\documentclass[10pt]{amsart}

\usepackage{amsmath,amsfonts,amssymb,amsthm,amscd,epsfig}

\def\C{\mathbb{C}}
\def\Z{\mathbb{Z}}

\def\g{\ensuremath{\mathfrak{g}}}
\def\n{\mathfrak{n}}

\def\V{\mathbf{V}}
\def\W{\mathbf{W}}
\def\v{\mathbf{v}}
\def\w{\mathbf{w}}
\def\e{\mathbf{e}}

\def\F{\mathcal{F}}

\def\L{\mathcal{L}}

\def\I{\mathbf{I}}
\def\i{\mathbf{i}}
\def\j{\mathbf{j}}
\def\a{\mathbf{a}}
\def\n{\mathbf{n}}
\def\ke{{\tilde e}}
\def\kf{{\tilde f}}

\DeclareMathOperator{\im}{Im} 
\DeclareMathOperator{\Hom}{Hom}

\DeclareMathOperator{\End}{End}
\DeclareMathOperator{\Aut}{Aut}
\DeclareMathOperator{\inc}{in}
\DeclareMathOperator{\out}{out}

\DeclareMathOperator{\tr}{tr}

\DeclareMathOperator{\wt}{wt}

\DeclareMathOperator{\Coker}{Coker}

\DeclareMathOperator{\diag}{diag}

\newtheorem{theo}{Theorem}[section]
\newtheorem{prop}[theo]{Proposition}
\newtheorem{lem}[theo]{Lemma}
\newtheorem{cor}[theo]{Corollary}

\newtheorem{defin}[theo]{Definition}
\newtheorem*{rem*}{Remark}
\newtheorem{rem}[theo]{Remark}
\numberwithin{equation}{section}

\begin{document}
\title[Geometric construction of crystal graphs]{A geometric construction
of crystal graphs using quiver varieties: extension to the
non-simply laced case}
\author{Alistair Savage}
\address{The Fields Institute for Research in Mathematical Sciences
and University of Toronto \\
Toronto, Ontario \\ Canada} \email{alistair.savage@aya.yale.edu}
\thanks{This research was supported by the Natural
Sciences and Engineering Research Council (NSERC) of Canada}
\date{June 3, 2004}

\begin{abstract}
We consider a generalization of the quiver varieties of Lusztig
and Nakajima to the case of all symmetrizable Kac-Moody Lie
algebras. To deal with the non-simply laced case one considers
admissible automorphisms of a quiver and the irreducible
components of the quiver varieties fixed by this automorphism.  We
define a crystal structure on these irreducible components and
show that the crystals obtained are isomorphic to those associated
to the crystal bases of the lower half of the universal enveloping
algebra and the irreducible highest weight representations of the
non-simply laced Kac-Moody Lie algebra.  As an application, we
realize the crystal of the spin representation of
$\mathfrak{so}_{2n+1}$ on the set of self-conjugate Young diagrams
that fit inside an $n \times n$ box.
\end{abstract}

\maketitle

\section*{Introduction}
In \cite{L91}, Lusztig associated various varieties to the quiver
(oriented graph) whose underlying graph is the Dynkin graph of a
simply laced Kac-Moody algebra (that is, a Kac-Moody algebra with
symmetric Cartan matrix).  Using these quiver varieties, he gave a
geometric realization of $U_q^-(\g)$, the lower half of the
quantum universal enveloping algebra of $\g$, and defined a
canonical basis \cite{L90} with remarkable properties.  Later,
Nakajima extended this quiver variety construction to yield a
geometric realization of the entire universal enveloping algebra
and its irreducible highest weight representations in the homology
of a new class of quiver varieties.  This geometric approach
allows one to prove such results as positivity and integrality of
the canonical basis that are difficult to prove algebraically.

In \cite{K91}, Kashiwara constructed the crystal base and the
global crystal base of $U_q^-(\g)$ and the highest weight
representations of $U_q(\g)$ in an algebraic way.  It was shown by
Grojnowski and Lusztig \cite{GL93} that the global crystal base
and the canonical basis coincide.  Kashiwara and Saito
\cite{KS97,S02} defined a crystal structure on the set of
irreducible components of Lusztig's and Nakajima's quiver
varieties and showed that the crystal obtained is isomorphic to
that associated to the crystal base of the lower half of the
universal enveloping algebra and the irreducible highest weight
representation respectively. There are other, purely
combinatorial, constructions of these crystals as well.  In
\cite{Sav03b}, the author described an explicit isomorphism
between the geometric construction from quiver varieties and the
combinatorial construction via Young tableaux and Young walls.

One limitation of the quiver variety approach to Kac-Moody Lie
algebras and their representations is that the geometric Lie
algebra action is only defined for those algebras whose
generalized Cartan matrix is symmetric. Thus, even for the finite
case, only types A, D, and E are covered.  However, the other
algebraic structures mentioned above, such as global bases and
crystal graphs, are defined for the more general case of
symmetrizable Kac-Moody Lie algebras.  In the current paper, we
consider an extension of the geometric construction of crystal
bases to this more general setting.  We consider quiver varieties
for symmetrizable Kac-Moody Lie algebras (which are not
necessarily simply laced) and introduce a crystal structure on
their irreducible components. In particular, we can treat all of
the finite dimensional simple Lie algebras. We give here a brief
overview of the construction.

Associated to a symmetrizable Kac-Moody Lie algebra (or its
generalized Cartan matrix), one can associate a valued graph
$\Gamma$. Let $\g(\Gamma)$ denote the Kac-Moody Lie algebra.  It
is known \cite[Prop 14.1.2]{L93} that any such valued graph can be
obtained from a quiver $Q$ along with an admissible automorphism
$\a$ of this graph. This admissible automorphism induces an
automorphism on the set of irreducible components of the quiver
varieties associated to $Q$. Let $\w$ be a dominant integral
weight of $\g(Q)$, the Kac-Moody Lie algebra with Dynkin diagram
the underlying graph of $Q$, which is fixed by the automorphism of
the weight lattice induced by $\a$. Then let $B_Q^g(\infty)$ and
$B_Q^g(\w)$ be the set of irreducible components of Lusztig's and
Nakajima's quiver varieties attached to the quiver $Q$ and let
$B_\Gamma^g(\infty)$ and $B_\Gamma^g(\w)$ be the subsets that are
fixed by the automorphism $\a$. Let $\I$ be the set of $\a$-orbits
of vertices of $Q$. It turns out that for $\i \in \I$ and $i,j \in
\i$, the Kashiwara operators $\ke_i$ and $\ke_j$ (resp. $\kf_i$
and $\kf_j$) which act on the sets $B_Q^g(\infty)$ and $B_Q^g(\w)$
commute and the operators $\ke_\i = \prod_{i \in \i} \ke_i$ and
$\kf_\i = \prod_{i \in \i} \kf_i$ act on the sets
$B_\Gamma^g(\infty)$ and $B_\Gamma^g(\w)$.  These operators (along
with the other maps required to define a crystal structure) endow
these sets with the structure of the crystal associated to the
crystal base of $U_q^-(\g(\Gamma))$ and $V_\Gamma(\w)$, the
irreducible highest weight representation of $U_q(\g(\Gamma))$ of
highest weight $\w$.  We can also define these operators directly,
in geometric terms, without reference to the operators $\ke_i$ and
$\kf_i$ for the quiver $Q$.

As an application of our construction, we describe an explicit
realization of the crystal graph of the spin representation of
$\mathfrak{so}_{2n+1}$ on the set of self-conjugate Young diagrams
that fit inside an $n \times n$ box.  Each Young diagram
corresponds to an irreducible component of the corresponding
quiver variety.  It turns out that we can actually extend this
crystal action to the full action of the Lie algebra and realize
the spin representation on the vector space spanned by these Young
diagrams.  This is similar to the realization of the spin
representations of $\mathfrak{so}_{2n}$ obtained in \cite{S03}.

The fact that the invariant irreducible components
$B_\Gamma^g(\w)$ of Nakajima's quiver varieties enumerate a basis
for highest weight representations of the corresponding
representation of the non-simply laced Kac-Moody Lie algebra was
demonstrated by F. Xu \cite{Xu97}. The proof in \cite{Xu97} uses
perverse sheaves while our method avoids their use and is almost
entirely crystal theoretic. The crystal structure on
$B_\Gamma^g(\infty)$ was mentioned in \cite{Xu97} but not defined.
A key step in our construction is the embedding of crystals of
non-simply laced type into those of simply laced type. We should
note that this has been done in some specific cases for Lie
algebras of affine type (see, for example, \cite{OSS03a,OSS03b}).
Our method works for arbitrary symmetrizable Kac-Moody algebras.

If we define the quiver variety associated to the valued graph
$\Gamma$ to be the union of the irreducible components of
Lusztig's (resp. Nakajima's) quiver variety associated to $Q$ that
are fixed by the isomorphism $\a$, then by the above we know that
the irreducible components of this new quiver variety are in 1-1
correspondence with a basis for $U_q^-(\g(\Gamma))$ (resp.
$V_\Gamma(\w)$).  Thus, the top dimensional Borel-Moore homology
of these quiver varieties has the same dimension (in each weight
space) as $U_q^-(\g(\Gamma))$ (resp. $V_\Gamma(\w)$). Therefore,
one should be able to define the action of $U_q(\g(\Gamma))$ on
these spaces via convolution in homology and obtain a geometric
realization of $U_q(\g(\Gamma))$ and $V_\Gamma(\w)$ as in
\cite{N98}.

One immediate application of the current work is an extension of
the results of \cite{Sav03b}.  The combinatorial methods of Young
tableaux and Young walls exist for non-simply laced Kac-Moody
algebras and thus the connection between the geometric and
combinatorial realizations of the crystal graphs should be able to
be extended to the more general setting. Ideally one should be
able to describe an explicit isomorphism between the geometric
realization and the path model of Littelmann \cite{Lit97} which
exists in the full generality of symmetrizable Kac-Moody Lie
algebras.

The organization of the paper is as follows.  In
Sections~\ref{sec:lus_def} and \ref{sec:def_nak} we review the
definitions of Lusztig's and Nakajima's quiver varieties and in
Sections~\ref{sec:lus-crystal} and \ref{sec:ca_qv} we recall the
crystal action on their irreducible components.  In
Section~\ref{sec:uev-aa} we prove a result about realizing the
crystal of $U_q^-(\g(\Gamma))$ inside the crystal of
$U_q^-(\g(Q))$.  We describe the geometric realization of the
crystal of $U_q^-(\g(\Gamma))$ in Section~\ref{sec:geom-real-lus}.
In Section~\ref{sec:rep-aa} we show how to realize the crystal of
$V_\Gamma(\w)$ inside that for $V_Q(\w)$ and in
Section~\ref{sec:geom-real-nak} we develop the geometric
realization of the crystal of $V_\Gamma(\w)$.  Finally, in
Section~\ref{sec:spinrep} we realize the crystal graph of the spin
representation of $\mathfrak{so}_{2n+1}$.

The author would like to thank I.B. Frenkel and C.M. Ringel for
useful discussions and suggestions.  He would also like to thank
H. Nakajima for bringing the reference \cite{Xu97} to his
attention.


\section{Lusztig's quiver variety}
\label{sec:lus_def}

In this section, we will recount the description given in
\cite{L91} of Lusztig's quiver varieties. See this reference for
details, including proofs.

Let $I$ and $E$ be the set of vertices and (unoriented) edges of
the Dynkin graph of a Kac-Moody Lie algebra with symmetric Cartan
matrix and let $H$ be the set of pairs consisting of an edge
together with an orientation of it. We denote the corresponding
quiver by $Q=(I,H)$ and the Kac-Moody Lie algebra by $\g(Q)$. For
$h \in H$, let $\inc(h)$ (resp. $\out(h)$) be the incoming (resp.
outgoing) vertex of $h$. We define the involution $\bar{\ }: H \to
H$ to be the function which takes $h \in H$ to the element of $H$
consisting of the same edge with opposite orientation.  An
\emph{orientation} of our graph is a choice of a subset $\Omega
\subset H$ such that $\Omega \cup \bar{\Omega} = H$ and $\Omega
\cap \bar{\Omega} = \emptyset$.

Let $\mathcal{V}$ be the category of finite-dimensional $I$-graded
vector spaces $\V = \oplus_{i \in I} \V_i$ over $\C$ with
morphisms being linear maps respecting the grading.  Then $\V \in
\mathcal{V}$ shall denote that $\V$ is an object of $\mathcal{V}$.
We identify the graded dimension $\v$ of $\V$ with the element
$\sum_{i \in I} \v_i \alpha_i$ of the root lattice of $\g(Q)$.
Here the $\alpha_i$ are the simple roots corresponding to the
vertices of our quiver.

Given $\V \in \mathcal{V}$, let
\[
\mathbf{E}(\v) = \bigoplus_{h \in H} \Hom (\V_{\out(h)},
\V_{\inc(h)}).
\]
Note that this space (and the others that follow) depend only on
the dimension of $\V$ up to isomorphism.  This is why we label the
space by the dimension $\v$.  For any subset $H' \subset H$, let
$\mathbf{E}_{H'}(\v)$ be the subspace of $\mathbf{E}(\v)$
consisting of all vectors $x = (x_h)$ such that $x_h=0$ whenever
$h \not\in H'$.  The algebraic group $G_\v = \prod_i \Aut(\V_i)$
acts on $\mathbf{E}(\v)$ and $\mathbf{E}_{H'}(\v)$ by
\[
(g,x) = ((g_i), (x_h)) \mapsto (g_{\inc(h)} x_h g_{\out(h)}^{-1}).
\]

Define the function $\varepsilon : H \to \{-1,1\}$ by $\varepsilon
(h) = 1$ for all $h \in \Omega$ and $\varepsilon(h) = -1$ for all
$h \in {\bar{\Omega}}$.  The Lie algebra of $G_\v$ is
$\mathbf{gl_v} = \prod_i \End(\V_i)$ and it acts on
$\mathbf{E}(\v)$ by
\[
(a,x) = ((a_i), (x_h)) \mapsto [a,x] = (x'_h) = (a_{\inc(h)}x_h -
x_h a_{\out(h)}).
\]
Let $\left<\cdot,\cdot\right>$ be the nondegenerate,
$G_\v$-invariant, symplectic form on $\mathbf{E}(\v)$ with values
in $\C$ defined by
\[
\left<x,y\right> = \sum_{h \in H} \varepsilon(h) \tr (x_h
y_{\bar{h}}).
\]
Note that $\mathbf{E}(\v)$ can be considered as the cotangent
space of $\mathbf{E}_\Omega(\v)$ under this form.

The moment map associated to the $G_{\v}$-action on the symplectic
vector space $\mathbf{E}(\v)$ is the map $\psi : \mathbf{E}(\v)
\to \mathbf{gl_v}$ with $i$-component $\psi_i : \mathbf{E}(\v) \to
\End \V_i$ given by
\[
\psi_i(x) = \sum_{h \in H,\, \inc(h)=i} \varepsilon(h) x_h
x_{\bar{h}} .
\]

\begin{defin}[\cite{L91}]
\label{def:nilpotent} An element $x \in \mathbf{E}(\v)$ is said to
be \emph{nilpotent} if there exists an $N \ge 1$ such that for any
sequence $h_1, h_2, \dots, h_N$ in $H$ satisfying $\out (h_1) =
\inc (h_2)$, $\out (h_2) = \inc (h_3)$, \dots, $\out (h_{N-1}) =
\inc (h_N)$, the composition $x_{h_1} x_{h_2} \dots x_{h_N} :
\V_{\out (h_N)} \to \V_{\inc (h_1)}$ is zero.
\end{defin}

\begin{defin}[\cite{L91}]
Let $\mathbf{E}^0(\v)$ be the set of elements $x \in
\mathbf{E}(\v)$ such that $\psi_i(x) = 0$ for all $i \in I$. Let
$\Lambda(\v)$ be the subset of $\mathbf{E}^0(\v)$ consisting of
nilpotent elements.
\end{defin}

\begin{prop}[\cite{L91}]
\label{prop:irrcomp-finite} For $\g(Q)$ a Kac-Moody Lie algebra of
finite type with symmetric Cartan matrix, the irreducible
components of $\Lambda(\v)$ are the closures of the conormal
bundles of the various $G_\v$-orbits in $\mathbf{E}_\Omega(\v)$.
\end{prop}


\section{Nakajima's quiver variety}
\label{sec:def_nak}

We introduce here a description of the quiver varieties first
presented in \cite{N94}.  See \cite{N94} and \cite{N98} for
details.

\begin{defin}[\cite{N94}]
\label{def:lambda} For $\v, \w \in (\Z_{\ge 0})^I$, choose
$I$-graded vector spaces $\V$ and $\W$ of graded dimensions $\v$
and $\w$ respectively.  We associate $\w$ with the element $\sum_i
\w_i \omega_i$ of the root lattice of $\g(Q)$, where the
$\omega_i$ are the fundamental weights of $\g(Q)$.  Recall that we
identified $\v$ with the weight $\sum_i \v_i \alpha_i$.  Then
define
\[
\Lambda \equiv \Lambda(\v;\w) = \Lambda(\v) \times \bigoplus_{i
\in I} \Hom (\V_i, \W_i).
\]
\end{defin}

Now, suppose that $\mathbf{S}$ is an $I$-graded subspace of $\V$.
For $x \in \Lambda(\v)$ we say that $\mathbf{S}$ is
\emph{$x$-stable} if $x(\mathbf{S}) \subset \mathbf{S}$.

\begin{defin}[\cite{N94}]
\label{def:lambda-stable} Let $\Lambda^{\text{st}} =
\Lambda(\v;\w)^{\text{st}}$ be the set of all $(x, t) \in
\Lambda(\v;\w)$ satisfying the following condition:  If
$\mathbf{S}=(\mathbf{S}_i)$ with $\mathbf{S}_i \subset \V_i$ is
$x$-stable and $t_i(\mathbf{S}_i) = 0$ for all $i \in I$, then
$\mathbf{S}_i = 0$ for all $i \in I$.
\end{defin}

The group $G_\v$ acts on $\Lambda(\v;\w)$ via
\[
(g,(x,t)) = ((g_i), ((x_h), (t_i))) \mapsto ((g_{\inc (h)} x_h
g_{\out (h)}^{-1}), (t_i g_i^{-1})).
\]
and the stabilizer of any point of $\Lambda(\v;\w)^{\text{st}}$ in
$G_{\v}$ is trivial (see \cite[Lemma~3.10]{N98}).  We then make
the following definition.
\begin{defin}[\cite{N94}]
\label{def:L} Let $\mathcal{L} \equiv \mathcal{L}(\v,\w) =
\Lambda(\v;\w)^{\text{st}} / G_{\v}$.
\end{defin}


\section{Crystal action on Lusztig's quiver varieties}
\label{sec:lus-crystal}

In this section we recall the crystal action on the set of
irreducible components of Lusztig's quiver varieties.
See~\cite{KS97} for details, including proofs.

Let $\v,{\bar \v},\v' \in (\Z_{\ge 0})^I$ such that $\v = {\bar
\v} + \v'$ and let $\V, {\bar \V}, \V'$ have dimensions $\v, {\bar
\v}, \v'$ respectively.  Consider the maps
\begin{equation} \label{eq:lus-proj1}
\mathbf{E}^0(\bar \v) \times \mathbf{E}^0(\v')
\stackrel{q_1}{\longleftarrow} \mathbf{E}^0({\bar \v},\v')
\stackrel{q_2}{\longrightarrow} \mathbf{E}^0(\v),
\end{equation}
where $\mathbf{E}^0({\bar \v},\v')$ is the variety of $(x,{\bar
\phi},\phi')$, where $x \in \mathbf{E}^0(\v)$ and ${\bar \phi} =
({\bar \phi}_i)$, $\phi' = (\phi'_i)$ give an exact sequence
\[
0 \longrightarrow {\bar \V}_i \stackrel{{\bar
\phi}_i}{\longrightarrow} \V_i \stackrel{\phi'_i}{\longrightarrow}
\V'_i \longrightarrow 0
\]
such that $\im {\bar \phi}$ is stable under $x$.  Thus $x$ induces
${\bar x} : {\bar \V} \to {\bar \V}$ and $x' : \V' \to \V'$.  The
maps $q_1$ and $q_2$ are defined by $q_1(x,{\bar \phi},\phi') =
(\bar x,x')$ and $q_2(x,{\bar \phi},\phi') = x$.

\begin{lem}[\cite{KS97}]
The following conditions are equivalent.
\begin{enumerate}
\item $x$ is nilpotent

\item Both $\bar x$ and $x'$ are nilpotent.
\end{enumerate}
\end{lem}

Thus \eqref{eq:lus-proj1} induces the maps
\begin{equation} \label{eq:lus-proj2}
\Lambda({\bar \v}) \times \Lambda(\v')
\stackrel{q_1}{\longleftarrow} \Lambda'({\bar \v},\v')
\stackrel{q_2}{\longrightarrow} \Lambda(\v),
\end{equation}
where $\Lambda'({\bar \v},\v') = q_2^{-1}(\Lambda(\v)) =
q_1^{-1}(\Lambda(\bar \v) \times \Lambda(\v'))$.

For $i \in I$ and $p \in \Z_{\ge 0}$, let
\[
\varepsilon_i(x) = \dim \Coker \left( \bigoplus_{h\, :\,
\inc(h)=i} V_{\out(h)} \stackrel{(x_h)}{\longrightarrow} V_i
\right),
\]
and
\[
\mathbf{E}^0(\v)_{i,p} = \{x \in \mathbf{E}^0(\v)\ |\
\varepsilon_i(x)=p\}.
\]
Then $\mathbf{E}^0(\v)_{i,p}$ is a locally closed subvariety of
$\mathbf{E}^0(\v)$.

Now assume that $\v = {\bar \v} + c\alpha_i$ for $c \in \Z_{\ge
0}$ and consider \eqref{eq:lus-proj1}.  One easily sees that
$\mathbf{E}^0(c\alpha_i) = \{0\}$.  Thus we have
\[
\mathbf{E}^0(\bar \v) \cong \mathbf{E}^0(\bar \v) \times
\mathbf{E}^0(c \alpha_i) \stackrel{\varpi_1}{\longleftarrow}
\mathbf{E}^0(\bar \v, c \alpha_i)
\stackrel{\varpi_2}{\longrightarrow} \mathbf{E}^0(\v).
\]
For $p \in \Z_{\ge 0}$ we have
\[
\varpi_1^{-1}(\mathbf{E}^0(\bar \v)_{i,p}) =
\varpi_2^{-1}(\mathbf{E}^0(\v)_{i,p+c}).
\]
Thus we define
\[
\mathbf{E}^0(\bar \v, c\alpha_i)_{i,p} =
\varpi_1^{-1}(\mathbf{E}^0(\bar \v)_{i,p}) = \varpi_2^{-1}
(\mathbf{E}^0(\v)_{i,p+c}).
\]
Setting $p=0$ we have the following diagram
\[
\mathbf{E}^0(\bar \v) \supset \mathbf{E}^0(\bar \v)_{i,0}
\stackrel{\varpi_1}{\longleftarrow} \mathbf{E}^0(\bar \v, c
\alpha_i)_{i,0} \stackrel{\varpi_2}{\longrightarrow}
\mathbf{E}^0(\v)_{i,c} \subset \mathbf{E}^0(\v).
\]
Note that $\mathbf{E}^0(\bar \v)_{i,0}$ is an open subvariety of
$\mathbf{E}^0(\bar \v)$.

\begin{lem}[\cite{KS97}]
We have the following.
\begin{enumerate}
\item $\varpi_2 : \mathbf{E}^0(\bar \v, c \alpha_i)_{i,0} \to
\mathbf{E}^0(\v)_{i,c}$ is a principal fiber bundle with $GL(\C^c)
\times \prod_{j \in I} GL({\bar \V}_j)$ as fiber.

\item $\varpi_1 : \mathbf{E}^0(\bar \v, c \alpha_i)_{i,0} \to
\mathbf{E}^0(\bar \v)_{i,0}$ is a smooth map whose fiber is a
connected rational variety.
\end{enumerate}
\end{lem}

Let $B^g_Q(\v,\infty)$ be the set of irreducible components of
$\Lambda(\v)$.  For $X \in B^g_Q(\v,\infty)$, we define
$\varepsilon_i(X) = \varepsilon_i(x)$ for a generic point $x$ of
$X$.  For $p \in \Z_{\ge 0}$, let $B^g_Q(\v,\infty)_{i,p}$ denote
the set of all elements $X$ of $B^g_Q(\v,\infty)$ such that
$\varepsilon_i(X)=p$.  From the above lemma, we obtain the
following.
\begin{prop}[\cite{KS97}]
We have
\[
B^g_Q(\bar \v;\infty)_{i,0} \cong B^g_Q(\v,\infty)_{i,c}.
\]
\end{prop}

Suppose that ${\bar X} \in B^g_Q(\bar \v;\infty)_{i,0}$
corresponds to $X \in B^g_Q(\v,\infty)_{i,c}$ by the above
isomorphism.  Then we define
\begin{gather*}
\kf_i^c : B^g_Q(\bar \v;\infty)_{i,0} \to
B^g_Q(\v,\infty)_{i,c},\quad \kf_i^c(\bar X) = X, \\
\ke_i^c : B^g_Q(\v,\infty)_{i,c} \to B^g_Q(\bar
\v;\infty)_{i,0},\quad \ke_i^c(X) = {\bar X}.
\end{gather*}
We then define maps
\begin{gather*}
\ke_i : \bigsqcup_\v B^g_Q(\v,\infty) \to \bigsqcup_\v
B^g_Q(\v,\infty) \sqcup \{0\}, \\
\kf_i : \bigsqcup_\v B^g_Q(\v,\infty) \to \bigsqcup_\v
B^g_Q(\v,\infty),
\end{gather*}
as follows.  For $c > 0$, we define
\[
\ke_i : B^g_Q(\v,\infty)_{i,c} \stackrel{\ke_i^c}{\longrightarrow}
B^g_Q(\bar \v;\infty)_{i,0}
\stackrel{\kf_i^{c-1}}{\longrightarrow} B^g_Q(\v -
\alpha_i;\infty)_{i,c-1},
\]
and $\ke_i(X)=0$ for $X \in B^g_Q(\v,\infty)_{i,0}$.  We define
\[
\kf_i : B^g_Q(\v,\infty)_{i,c} \stackrel{\ke_i^c}{\longrightarrow}
B^g_Q(\infty;\bar \v)_{i,0}
\stackrel{\kf_i^{c+1}}{\longrightarrow} B^g_Q(\infty; \v +
\alpha_i)_{i,c+1}.
\]
Furthermore, we define a map
\[
\wt : \bigsqcup_\v B^g_Q(\v,\infty) \to P,\quad \wt(X) = - \sum_{i
\in I} \v_i \alpha_i \text{ for } X \in B^g_Q(\v,\infty)
\]
and we set
\[
\varphi_i(X) = \varepsilon_i(X) + \left< h_i, \wt(X) \right>.
\]
Let $B^g_Q(\infty) = \bigsqcup_\v B^g_Q(\v,\infty)$.

\begin{prop}[\cite{KS97}]
The maps defined above make $B^g_Q(\infty)$ a crystal and it is
isomorphic to the crystal $B_Q(\infty)$ associated to the crystal
base of $U_q^-(\g(Q))$.
\end{prop}


\section{Crystal action on Nakajima's quiver varieties}
\label{sec:ca_qv} In this section, we review the realization of
the crystal graph of integrable highest weight representations of
$\g(Q)$ via quiver varieties. See \cite{S02} for details,
including proofs.

Let $\mathbf{w, v, v', v''} \in (\Z_{\ge 0})^I$ be such that $\v =
\mathbf{v'} + \mathbf{v''}$.  Consider the maps
\begin{equation}
\label{eq:diag_action} \Lambda(\v'';\mathbf{0}) \times
\Lambda(\v';\w) \stackrel{p_1}{\leftarrow} \mathbf{\tilde F
(v,w;v'')} \stackrel{p_2}{\rightarrow} \mathbf{F(v,w;v'')}
\stackrel{p_3}{\rightarrow} \Lambda(\v;\w),
\end{equation}
where the notation is as follows.  A point of
$\mathbf{F(v,w;v'')}$ is a point $(x,t) \in \Lambda(\v;\w)$
together with an $I$-graded, $x$-stable subspace $\mathbf{S}$ of
$\V$ such that $\dim \mathbf{S} = \mathbf{v'} = \v - \v''$.  A
point of $\mathbf{\tilde
  F (v,w;v'')}$ is a point $(x,t,\mathbf{S})$ of $\mathbf{F(v,w;v'')}$
together with a collection of isomorphisms $R'_i : \V'_i \cong
\mathbf{S}_i$ and $R''_i : \V''_i \cong \V_i / \mathbf{S}_i$ for
each $i \in I$.  Then we define $p_2(x,t,\mathbf{S}, R',R'') =
(x,t,\mathbf{S})$, $p_3(x,t,\mathbf{S}) = (x,t)$ and
$p_1(x,t,\mathbf{S},R',R'') = (x'',x',t')$ where $x'', x', t'$ are
determined by
\begin{align*}
R'_{\inc(h)} x'_h &= x_h R'_{\out(h)} : \V'_{\out(h)} \to
\mathbf{S}_{\inc(h)}, \\
t'_i &= t_i R'_i : \V'_i \to \W_i \\
R''_{\inc(h)} x''_h &= x_h R''_{\out(h)} : \V''_{\out(h)} \to
\V_{\inc(h)} / \mathbf{S}_{\inc(h)}.
\end{align*}
It follows that $x'$ and $x''$ are nilpotent.

\begin{lem}[{\cite[Lemma 10.3]{N94}}]
One has
\[
(p_3 \circ p_2)^{-1} (\Lambda(\v;\w)^{\text{st}}) \subset p_1^{-1}
(\Lambda(\v'';\mathbf{0}) \times \Lambda(\v';\w)^{\text{st}}).
\]
\end{lem}

Thus, we can restrict \eqref{eq:diag_action} to
$\Lambda^{\text{st}}$, forget the
$\Lambda(\v'';\mathbf{0})$-factor and consider the quotient by
$G_\v$, $G_\mathbf{v'}$.  This yields the diagram
\begin{equation}
\label{eq:diag_action_mod} \mathcal{L}(\v', \w)
\stackrel{\pi_1}{\leftarrow} \mathcal{F}(\v, \w; \v - \mathbf{v'})
\stackrel{\pi_2}{\rightarrow} \mathcal{L}(\v, \w),
\end{equation}
where
\[
\mathcal{F}(\v, \w; \v - \mathbf{v'}) \stackrel{\text{def}}{=} \{
(x,t,\mathbf{S}) \in \mathbf{F(\v,\w;\v-\v')}\,
  |\, (x,t) \in \Lambda(\v;\w)^{\text{st}} \} / G_\v.
\]

For $i \in I$ define $\varepsilon_i : \Lambda(\v; \w) \to \Z_{\ge
0}$ by
\[
\varepsilon_i((x,t)) = \dim_\C \Coker \left( \bigoplus_{h\, :\,
\inc(h)=i} V_{\out(h)} \stackrel{(x_h)}{\longrightarrow} V_i
\right).
\]
Then, for $c \in \Z_{\ge 0}$, let
\[
\mathcal{L}(\v,\w)_{i,c} = \{[x,t] \in \mathcal{L}(\v,\w)\ |\
\varepsilon_i((x,t)) = c\}
\]
where $[x,t]$ denotes the $G_\v$-orbit through the point $(x,t)$.
We see that $\mathcal{L}(\v,\w)_{i,c}$ is a locally closed
subvariety of $\mathcal{L}(\v,\w)$.

Assume $\mathcal{L}(\v,\w)_{i,c} \ne \emptyset$ and let $\v' = \v
- c\mathbf{e}^i$ where $\mathbf{e}^i_j = \delta_{ij}$.  Then
\[
\pi_1^{-1}(\mathcal{L}(\v',\w)_{i,0}) =
\pi_2^{-1}(\mathcal{L}(\v,\w)_{i,c}).
\]
Let
\[
\mathcal{F}(\v,\w;c\mathbf{e}^i)_{i,0} =
\pi_1^{-1}(\mathcal{L}(\v',\w)_{i,0}) =
\pi_2^{-1}(\mathcal{L}(\v,\w)_{i,c}).
\]
We then have the following diagram.
\begin{equation}
\label{eq:crystal-action} \mathcal{L}(\v',\w)_{i,0}
\stackrel{\pi_1}{\longleftarrow} \mathcal{F}(\v,\w;
c\mathbf{e}^i)_{i,0} \stackrel{\pi_2}{\longrightarrow}
\mathcal{L}(\v,\w)_{i,c}
\end{equation}
The restriction of $\pi_2$ to $\mathcal{F}(\v,\w;
c\mathbf{e}^i)_{i,0}$ is an isomorphism since the only possible
choice for the subspace $\mathbf{S}$ of $\V$ is to have
$\mathbf{S}_j = \mathbf{V}_j$ for $j \ne i$ and $\mathbf{S}_i$
equal to the sum of the images of the $x_h$ with $\inc(h)=i$.
$\mathcal{L}(\v',\w)_{i,0}$ is an open subvariety of
$\mathcal{L}(\v',\w)$.

\begin{lem}[\cite{S02}]
\begin{enumerate}
\item For any $i \in I$,
\[
\mathcal{L}(\mathbf{0},\w)_{i,c} =
\begin{cases}
pt & \text{if $c=0$} \\
\emptyset & \text{if $c >0$}
\end{cases}.
\]
\item Suppose $\mathcal{L}(\v,\w)_{i,c} \ne \emptyset$ and $\v' =
\v - c\mathbf{e}^i$.  Then the fiber of the restriction of $\pi_1$
to $\mathcal{F}(\v, \w; c\mathbf{e}^i)_{i,0}$ is isomorphic to a
Grassmanian variety.
\end{enumerate}
\end{lem}

\begin{cor}
\label{cor:irrcomp-isom} Suppose $\mathcal{L}(\v,\w)_{i,c} \ne
\emptyset$.  Then there is a 1-1 correspondence between the set of
irreducible components of $\mathcal{L}(\v - c\mathbf{e}^i,
\w)_{i,0}$ and the set of irreducible components of
$\mathcal{L}(\v, \w)_{i,c}$.
\end{cor}

Let $B^g_Q(\v,\w)$ denote the set of irreducible components of
$\mathcal{L}(\v,\w)$ and let $B^g_Q(\w) = \bigsqcup_\v
B^g_Q(\v,\w)$. For $X \in B^g_Q(\v,\w)$, let $\varepsilon_i(X) =
\varepsilon_i((x,t))$ for a generic point $[x,t] \in X$.  Then for
$p \in \Z_{\ge 0}$ define
\[
B^g_Q(\v,\w)_{i,p} = \{X \in B^g_Q(\v,\w)\ |\ \varepsilon_i(X) =
p\}.
\]
Then by Corollary~\ref{cor:irrcomp-isom}, $B^g_Q(\v -
c\mathbf{e}^i,\w)_{i,0} \cong B^g_Q(\v, \w)_{i,c}$.

Suppose that ${\bar X} \in B^g_Q(\v - c\mathbf{e}^i,\w)_{i,0}$
corresponds to $X \in B^g_Q(\v,\w)_{i,c}$ by the above
isomorphism. Then we define maps
\begin{gather*}
\kf_i^c : B^g_Q(\v - c\mathbf{e}^i,\w)_{i,0} \to
B^g_Q(\v,\w)_{i,c},\quad \kf_i^c({\bar X})
= X, \\
\ke_i^c : B^g_Q(\v,\w)_{i,c} \to B^g_Q(\v -
c\mathbf{e}^i,\w)_{i,0},\quad \ke_i^c(X) = {\bar X}.
\end{gather*}
We then define the maps
\[
\ke_i, \kf_i : B^g_Q(\w) \to B^g_Q(\w) \sqcup \{0\}
\]
by
\begin{gather*}
\ke_i : B^g_Q(\v,\w)_{i,c} \stackrel{\ke_i^c}{\longrightarrow}
B^g_Q(\v - c\mathbf{e}^i, \w)_{i,0}
\stackrel{\kf_i^{c-1}}{\longrightarrow} B^g_Q(\v -
\mathbf{e}^i, \w)_{i,c-1}, \\
\kf_i : B^g_Q(\v,\w)_{i,c} \stackrel{\ke_i^c}{\longrightarrow}
B^g_Q(\v - c\mathbf{e}^i, \w)_{i,0}
\stackrel{\kf_i^{c+1}}{\longrightarrow} B^g_Q(\v + \mathbf{e}^i,
\w)_{i,c+1}.
\end{gather*}
We set $\ke_i(X)=0$ for $X \in B^g_Q(\v,\w)_{i,0}$ and
$\kf_i(X)=0$ for $X \in B^g_Q(\v,\w)_{i,c}$ with
$B^g_Q(\v+\e^i,\w)_{i,c+1} = \emptyset$. We also define
\begin{gather*}
\wt : B^g_Q(\w) \to P,\quad \wt(X) = \sum_{i\in I} \left(
\mathbf{w}_i \omega_i - \mathbf{v}_i \alpha_i \right) \text{
  for } X \in B^g_Q(\v,\w), \\
\varphi_i(X) = \varepsilon_i(X) + \left< h_i, \wt(X) \right>.
\end{gather*}

Recall that we can consider $\w$ to be an dominant integral weight
by $\mathbf{w} = \sum_i \w_i \omega_i$.

\begin{prop}[\cite{S02}]
Under the maps defined above, $B^g_Q(\w)$ is a crystal and is
isomorphic to the crystal $B_Q(\w)$ associated to the crystal base
of $V_Q(\w)$, the highest weight $U_q(\g)$-module with highest
weight $\w$.
\end{prop}


\section{The crystal of the universal enveloping algebra and admissible automorphisms}
\label{sec:uev-aa}

As before, let $Q$ be a quiver without vertex loops.  The
corresponding symmetric generalized Cartan matrix is the matrix
$A$ indexed by $I$ with entries
\[
a_{ij} = \begin{cases} 2 & i=j \\
-\#\{\text{edges joining vertices $i$ and $j$}\} & i \ne j
\end{cases}.
\]
Note that by the word \emph{edge} we mean unoriented edge (i.e.
elements of $E$).  Again, we let $\g(Q)$ denote the associated
symmetric Kac-Moody algebra, with root system $\Delta(Q)$ (see
\cite{K}).

An admissible automorphism $\a$ of a (double) quiver $Q$ is an
automorphism of the underlying graph such that no edge connects
two vertices in the same $\a$-orbit. Following \cite{L93} we
construct a symmetric matrix $M$ indexed by the vertex $\a$-orbits
$\I$.  We let the $(\i,\j)$ entry be
\[
m_{\i\j} = \begin{cases} 2\#\{\text{vertices in $\i$th orbit}\} &
\i = \j \\
-\#\{\text{edges joining a vertex in $\i$th orbit and a vertex in
$\j$th orbit}\} & \i \ne \j \end{cases}.
\]
Then let
\[
d_\i = m_{\i \i}/2 = \#\{\text{vertices in $\i$th orbit}\}
\]
and set $D = \diag(d_\i)$.  Then $C = D^{-1}M$ is a symmetrizable
generalized Cartan matrix.  Let $\Gamma$ denote the corresponding
valued graph.  That is, $\Gamma$ has vertex set $\I$ and whenever
$c_{\i \j} \ne 0$, we draw an edge connecting $\i$ and $\j$
equipped with the ordered pair $(|c_{\j \i}|,|c_{\i \j}|)$.  It is
known \cite[Prop 14.1.2]{L93} that any symmetrizable generalized
Cartan matrix (and corresponding valued graph) can be obtained
from a pair $(Q,\a)$ in this way.  Note that the fact that $\a$ is
admissible ensures that $\Gamma$ has no vertex loops.  Let
$\g(\Gamma)$ be the Kac-Moody algebra associated to $C$, with root
system $\Delta(\Gamma)$.

Let $(-,-)_Q$ and $(-,-)_\Gamma$ be the symmetric bilinear forms
determined by the matrices $A$ and $M$ respectively.  The
automorphism $\a$ acts naturally on the root lattice $\Z I$ for
$Q$, and $(-,-)_Q$ is $\a$-invariant.  There is a canonical
bijection
\[
f: (\Z I)^\a \to \Z \I,\quad f(\alpha)_\i = \alpha_i \text{ for
any } i \in \i,
\]
from the fixed points in the root lattice for $Q$ to the root
lattice for $\Gamma$.  We will often suppress the bijection $f$
and consider the root lattice of $\Gamma$ to be the fixed points
in the root lattice for $Q$.  In particular, we have the simple
roots for $\Gamma$ given by
\[
\alpha_\i = \sum_{i \in \i} \alpha_i.
\]
We also define
\[
h_\i = \frac{1}{d_\i} \sum_{i \in \i} h_i.
\]
Then the entries of $C$ are given by $c_{\i \j} = \left<
\alpha_\i, h_\j \right>$.

Recall that $B_Q(\infty)$ and $B_Q^g(\infty)$ are the algebraic
and geometric crystals of $U_q^-(\g(Q))$ respectively. We know
that $B_Q(\infty) \cong B_Q^g(\infty)$. One of the goals of this
paper is to define a geometric crystal $B_\Gamma^g(\infty)$ of
$U_q^-(\g(\Gamma))$.

We use the same notation $\wt$, $\varepsilon_i$, $\varphi_i$,
$\ke_i$ and $\kf_i$, $i \in I$, for the maps of the crystals
$B_Q(\infty)$ and $B_Q^g(\infty)$.

\begin{lem}
\label{lem:ef-comm}
Suppose $i$ and $j$ are in the same vertex
orbit $\i$. Then
\[
\ke_i \ke_j = \ke_j \ke_i,\quad \kf_i \kf_j = \kf_j \kf_i.
\]
\end{lem}

\begin{proof}
We first briefly recall the definition of the Kashiwara operators
$\ke_i$ and $\kf_i$ in the algebraic setting.  For any $P \in
U_q^-(\g(Q))$, there exist unique $Q,R \in U_q^-(\g(Q))$ such that
\[
[e_i,P] = \frac{q^{h_i} Q - q^{-h_i}R}{q - q^{-1}}.
\]
We can therefore define an endomorphism $e_i'$ of $U_q^-(\g(Q))$
given by $e_i'(P) = R$. According to \cite{K91} we have
\[
U_q^-(\g(Q)) = \bigoplus_{n \ge 0} f_i^{(n)} \ker e_i'.
\]
Recall that $f_i^{(n)} = f_i^n/[n]_i!$ where $[n]_i = (q^n -
q^{-n})/(q - q^{-1})$ and $[n]_i! = \prod_{k=1}^n [k]_i$. Then the
Kashiwara operators are induced by the operators on $U_q^-(\g(Q))$
given by
\begin{gather*}
\kf_i (f_i^{(n)} u) = f_i^{(n+1)} u, \\
\ke_i (f_i^{(n)} u) = f_i^{(n-1)} u,
\end{gather*}
for $u \in \ker e_i'$, where $f_i^{(n)} = 0$ for $n < 0$.

Fix $w \in U_q^-(\g(Q))$ and let $w = f_i^{(n)} u$ for $u \in \ker
e_i'$.  Let $u = f_j^{(m)} v$ for $v \in \ker e_j'$.  Since $i$
and $j$ lie in the same $\a$-orbit and $\a$ is an admissible
automorphism, there is no edge joining the vertices $i$ and $j$.
Therefore $[f_i,f_j] = [e_i,f_j] = [e_j,f_i]=0$ and so we have
\[
w = f_i^{(n)} f_j^{(m)} v = f_j^{(m)} f_i^{(n)} v, \quad v \in
\ker e_i' \cap \ker e_j'.
\]
Also, since $[e_i,f_j] = [e_j,f_i]=0$, we have $f_i^{(l)} v \in
\ker e_j'$ and $f_j^{(l)} v \in \ker e_i'$ for all $l$.
Therefore,
\begin{align*}
\ke_i \ke_j w &= \ke_i \ke_j (f_j^{(m)} f_i^{(n)} v) \\
&= \ke_i (f_j^{(m-1)} f_i^{(n)} v) \\
&= \ke_i (f_i^{(n)} f_j^{(m-1)} v) \\
&= f_i^{(n-1)} f_j^{(m-1)} v \\
&= f_j^{(m-1)} f_i^{(n-1)} v \\
&= \ke_j (f_j^{(m)} f_i^{(n-1)} v) \\
&= \ke_j (f_i^{(n-1)} f_j^{(m)} v) \\
&= \ke_j \ke_i (f_i^{(n)} f_j^{(m)} v) \\
&= \ke_j \ke_i w.
\end{align*}
The proof that $\kf_i \kf_j = \kf_j \kf_i$ is analogous.
\end{proof}

By Lemma~\ref{lem:ef-comm}, we can unambiguously define the
operators
\[
\ke_\i = \prod_{i \in \i} \ke_i,\quad \kf_\i = \prod_{i \in \i}
\kf_i,\quad \i \in \I.
\]
For $b \in B_Q(\infty)$ and $\i \in \I$, we also define
\begin{gather*}
\varepsilon_\i (b) = \max \{k \ge 0\ |\ \ke_\i^k b \ne 0\}, \\
\varphi_\i (b) = \varepsilon_\i(b) + \left< h_\i, \wt(b)\right>.
\end{gather*}

Let $B_\Gamma(\infty)$ be the subset of $B_Q(\infty)$ generated by
the the $\kf_\i$, $\i \in \I$ acting on the highest weight element
$b_\infty$ of $B_Q(\infty)$.  Note that if we restrict the map
$\wt : B_Q(\infty) \to P(Q)$, where $P(Q)$ is the weight lattice
of $\g(Q)$, to the subset $B_\Gamma(\infty)$, then the image lies
in the subset of $P(Q)$ that is invariant under the natural action
of $\a$.  Thus we can view it as a map $\wt : B_\Gamma(\infty) \to
P(\Gamma)$ where $P(\Gamma)$ is the weight lattice of
$\g(\Gamma)$.

\begin{prop}
The set $B_\Gamma(\infty)$, along with the maps $\wt$,
$\varepsilon_\i$, $\varphi_\i$, $\ke_\i$ and $\kf_\i$ is a
$\g(\Gamma)$-crystal.
\end{prop}

\begin{proof}
We have to show that Axioms (1)-(7) of \cite[Definition 4.5.1]{HK}
defining a crystal are satisfied.  Axioms (1)-(6) follow easily
from the axioms for the crystal $B_Q(\infty)$ and Axiom (7) is
vacuous since we never have $\varphi_\i (b) = -\infty$.
\end{proof}

We want to show that the crystal $B_\Gamma(\infty)$ is actually
isomorphic to the crystal of $U_q^-(\g(\Gamma))$.

\begin{defin}
For $i \in I$, $B_i$ is the crystal defined as follows.
\begin{gather*}
B_i = \{b_i(n)\ |\ n \in \Z\}, \\
\wt(b_i(n)) = n\alpha_i, \\
\varphi_i(b_i(n)) = n,\quad \varepsilon_i(b_i(n)) = -n, \\
\varphi_j(b_i(n)) = \varepsilon_j(b_i(n)) = -\infty \text{ for } i
\ne j.
\end{gather*}
The action of $\ke_j$ and $\kf_j$ is given by
\begin{gather*}
\ke_i(b_i(n)) = b_i(n+1), \\
\kf_i(b_i(n)) = b_i(n-1), \\
\ke_j(b_i(n)) = \kf_j(b_i(n)) = 0,\ i \ne j.
\end{gather*}
Let $b_i = b_i(0)$.  We similarly define $B_\i$ for $\i \in \I$.
\end{defin}

We now consider an equivalent definition of $B_\i$ for $\i \in
\I$. Let $\i = \{i_1, \dots, i_k\}$ and define
\[
B'_\i = \{b'_\i (n) := b_{i_1}(n) \otimes \dots \otimes
b_{i_k}(n)\ |\ n \in \Z\} \subset B_{i_1} \otimes \dots \otimes
B_{i_k}.
\]
By the definition of the tensor product of crystals (see
\cite{HK}), we have that
\begin{gather*}
\varphi_i (b'_\i(n)) = n,\quad \varepsilon_i (b'_\i(n)) = -n
\text{
for } i \in \i \\
\varphi_j (b'_\i(n)) = \varepsilon_j (b'_\i (n)) = -\infty \text{
for } j \not \in \i \\.
\end{gather*}
Thus we define $\varepsilon_\i = \varepsilon_i$ and $\varphi_\i =
\varphi_i$ for any $i \in \i$.  Note also that
\[
\wt(b'_\i(n)) = \sum_{i \in \i} \wt(b_i(n)) = \sum_{i \in \i}
n\alpha_i = n \alpha_\i.
\]
It is also straightforward to verify from the definition of the
tensor product of crystals that
\begin{gather*}
\ke_\i (b'_\i(n)) = b'_\i (n+1) \\
\kf_\i (b'_\i(n)) = b'_\i (n-1) \\
\ke_\j (b'_\i(n)) = \kf_\j (b'_\i(n)) = 0,\quad \i \ne \j.
\end{gather*}
Thus we see that $B'_\i \cong B_\i$ and henceforth we identify the
two.  In particular, we identify $b'_\i(n)$ and $b_\i(n)$.

The following characterization of the crystal of a universal
enveloping algebra will be useful.

\begin{prop}[\cite{KS97}] \label{prop:crystal-char}
Let $B$ be a $\g(\Gamma)$-crystal and $b_0$ an element of $B$ with
weight zero.  Assume the following conditions.
\begin{enumerate}
\item \label{char1} $\wt(B) \subset Q(\Gamma)^-$

\item \label{char2} $b_0$ is the unique element of $B$ with weight
zero.

\item \label{char3} $\varepsilon_\i (b_0) = 0$ for every $\i$.

\item \label{char4} $\varepsilon_\i (b) \in \Z$ for any $b$ and
$\i$.

\item \label{char5} For every $\i$, there exists a strict
embedding $\Psi_\i : B \to B \otimes B_\i$.

\item \label{char6} $\Psi_\i(B) \subset B \times \{\kf_\i^n b_\i\
|\ n \ge 0\}$.

\item \label{char7} For any $b \in B$ such that $b \ne b_0$ there
exists $\i$ such that $\Psi_\i(b) = b' \otimes \kf_i^n b_\i$ with
$n >0$.
\end{enumerate}
Then $B$ is isomorphic to the crystal associated to the crystal
base of $U_q^-(\g(\Gamma))$.  The same result holds for
$U_q^-(\g(Q))$ if we replace $\I$ by $I$ (and $\i$ by $i$) in the
above.
\end{prop}

\begin{prop}
$B_\Gamma(\infty)$ is isomorphic to the crystal associated to the
crystal base of $U_q^-(\g(\Gamma))$.
\end{prop}

\begin{proof}
We prove this by an appeal to Proposition~\ref{prop:crystal-char}.
Let $b_0 = b_\infty$ be the unique element of $B_\Gamma(\infty)$
of weight zero. Conditions~\eqref{char1}-\eqref{char4} are
obvious. For $\i = \{i_1,\dots,i_d\} \in \I$, we define
\[
\Psi_\i = \Psi_{i_1} \circ \Psi_{i_2} \circ \dots \circ \Psi_{i_d}
: B \to B \otimes B_{i_1} \otimes B_{i_2} \otimes \dots \otimes
B_{i_d}.
\]
Now, by conditions~\eqref{char1}, \eqref{char2} and \eqref{char6}
for the case of $U_g^-(\g(Q))$, we see that $\Psi_i(b_\infty) =
b_\infty \otimes b_i$.  Thus
\[
\Psi_\i(b_0) = \Psi_\i(b_\infty) = b_\infty \otimes b_{i_1}
\otimes \dots \otimes b_{i_d} = b_\infty \otimes b_\i \in
B_\Gamma(\infty) \otimes B_\i.
\]
Thus, since each $\Psi_{i_k}$ is strict (that is, it commutes with
all $\ke_j$ and $\kf_j$), we see that the image of $\Psi_\i$
restricted to $B_\Gamma(\infty)$ lies in $B_\Gamma(\infty) \otimes
B_\i$.  We also denote the restriction by $\Psi_\i$.  Thus,
condition~\eqref{char5} is satisfied.

Condition~\eqref{char6} for $\Psi_\i$ now follows from the
corresponding conditions for the $\Psi_{i_k}$.  For $b \in
B_\Gamma(\infty)$ such that $b \ne b_0$, pick an $i$ such that
condition~\eqref{char7} holds for $U_q^-(\g(Q))$.  Then it follows
that condition~\eqref{char7} holds for $U_q^-(\g(\Gamma))$ for the
orbit $\i$ containing $i$.
\end{proof}


\section{A geometric realization of the crystal of the universal
enveloping algebra in the non-simply laced case}
\label{sec:geom-real-lus}

Since we know that the geometric crystal $B^g_Q(\infty)$ is
isomorphic to the algebraic crystal $B_Q(\infty)$, we can define
$B^g_\Gamma(\infty)$ just as we defined $B_\Gamma(\infty)$ and we
know that $B^g_\Gamma(\infty)$ is isomorphic to the crystal
associated to the crystal base of $U_q^-(\g(\Gamma))$.  By
definition, the underlying set of $B^g_\Gamma(\infty)$ consists of
the irreducible components of Lusztig's quiver variety generated
from the highest weight component (which consists of a point) by
the action of the operators $\kf_\i$ for $\i \in \I$.  We now seek
to describe this set in a more direct and geometric way.

As described in \cite{H04}, the automorphism $\a$ gives rise to an
autoequivalence of the category of representations of the quiver
as follows.  Recall our quiver $Q$ with vertices $I$, edges $E$,
and oriented edges $H$.  Note that the automorphism $\a$ acts
naturally on $H$. By a representation of the quiver $Q$, we mean a
collection of vector spaces $\V=(\V_i)_{i \in I}$ and linear maps
$(x_h : \V_{\out(h)} \to \V_{\inc(h)})_{h \in H}$.  The dimension
vector $\v = \dim \V \in \Z I$ is defined to be $\sum_{i \in I}
(\dim \V_i) \alpha_i$ where $\{\alpha_i\}_{i \in I}$ is viewed as
the standard basis of $\Z I$ or the set of simple roots depending
on the context.

Given a representation of $Q$, we define a new representation as
follows.  The new collection of vector spaces is $^\a \V$ where
$(^\a \V)_i = \V_{\a^{-1}(i)}$ and the maps are given by $(^\a
x)_h = x_{\a^{-1}(h)}$.  We also have a natural action on
morphisms between representations but we will not need this
action.  Thus we have a functor $F(\a)$ on the category of
representations.  The functor $F(\a)$ is an autoequivalence of
this category. Note that in \cite{H04}, the case of the single
quiver (one orientation per edge) is considered while we are
dealing with the double quiver here.

Note that $\dim {^\a \V} = \a(\dim \V)$.  So if we restrict
ourselves to those $\V$ whose dimensions are invariant under $\a$,
we can consider $F(\a)$ as an automorphism of the space of
representations.  This induces an automorphism of Lusztig's quiver
variety and thus an automorphism of the set of irreducible
components (we assume that our chosen orientation $\Omega$ is
compatible with the automorphism of the quiver which is always
possible by \cite[\S 12.1.1]{L93}). We denote the automorphism of
the set of irreducible components thus obtained by $\a$.

\begin{prop} \label{prop:univ-irrcomp-inv}
The underlying set of $B^g_\Gamma(\infty)$ is precisely the subset
of irreducible components of $B^g_Q(\infty)$ consisting of those
components that are invariant under the automorphism $\a$.
\end{prop}

One should note that the irreducible components invariant under
$\a$ are not necessarily invariant under $F(\a)$.  That is, each
individual point in the irreducible component is not necessarily
invariant under $F(\a)$.

\begin{proof}
Let $B'$ be the set of irreducible components invariant under
$\a$. It is easy to see that for $i \in I$, we have
\[
a \kf_i a^{-1} = \kf_{a(i)}.
\]
Thus for $\i \in \I$,
\begin{align*}
a \kf_\i a^{-1} &= a \left( \prod_{i \in \i} \kf_i \right) a^{-1}
= \prod_{i \in \i} (a \kf_i a^{-1}) \\
&= \prod_{i \in \i} \kf_{a(i)} = \prod_{i \in \i} \kf_i = \kf_\i.
\end{align*}
Therefore, since the highest weight element of
$B^g_\Gamma(\infty)$, which we will denote by $X_\infty$, is
invariant under $\a$, we have that $B^g_\Gamma(\infty) \subset B'$
(as sets).

Now, let $X$ be an element of $B^g_Q(\infty)$ such that $\a(X)=X$.
Let us order the set $\I$.  So we have $\I = \{\i_1,\dots,\i_l\}$
and we let $\i_p = \{i_{p1}, \dots, i_{p d_p}\}$.  Where $d_p =
d_{\i_p}$ is the number of vertices in the orbit $\i_p$. We define
a collection of non-negative integers $\n = (n_j^{b,c})_{1 \le j
\le k, 1 \le b \le l, 1 \le c \le d_b}$ inductively as follows.
\begin{align*}
n_1^{1,1} &= \varepsilon_{i_{11}}(X), \\
X_j^{b,c} &= \ke_{i_{bc}}^{n_j^{b,c}} \cdots
\ke_{i_{b1}}^{n_j^{b,1}} \cdots \ke_{i_{1d_1}}^{n_j^{1,d_1}}
\cdots \ke_{i_{11}}^{n_j^{1,1}} \cdots
\ke_{i_{ld_l}}^{n_1^{l,d_l}} \cdots \ke_{i_{l1}}^{n_1^{l,1}}
\cdots \ke_{i_{1d_1}}^{n_1^{1,d_1}} \cdots
\ke_{i_{11}}^{n_1^{1,1}} X, \\
n_j^{b,c} &= \varepsilon_{i_{bc}}(X_j^{b,c-1}) \text{ if } c>1, \\
n_j^{b,1} &= \varepsilon_{i_{b1}}(X_j^{b-1,d_{b-1}}) \text{ if } b>1, \\
n_j^{1,1} &= \varepsilon_{i_{11}}(X_{j-1}^{l,d_l}) \text{ if }
j>1.
\end{align*}
and $k$ is chosen such that $X_k^{l,d_l} = X_\infty$.  Intuitively
speaking, we are applying $\ke_i$ for $i \in \i_1$ as many times
as possible (without getting 0).  Then we apply $\ke_i$ for $i \in
\i_2$ as many times as possible, etc.  Once we have applied
$\ke_i$ for $i \in \i_l$, we start again at $\i_1$.  We continue
this process until we obtain $X_\infty$, which we will always do
since $B^g_Q(\infty)$ is connected.  Then the $n_j^{b,c}$ tell us
how many times we applied the various $\ke_i$.  Note that it does
not really matter what order we pick the various $i$ in each $\i$
since such $\ke_i$ commute.  Also note that $X$ uniquely
determines $\n$ and vice versa since using $\n$ we can apply the
operators $\kf_i$ to $X_\infty$ to obtain $X$.

Identify $\n$ with the sequence
\begin{align*}
\n = ({i_{l d_l}}^{n_k^{l,d_l}}, &\dots, {i_{l1}}^{n_k^{l,1}},
\dots, {i_{1 d_1}}^{n_k^{1,d_1}},\dots, {i_{11}}^{n_k^{1,1}},
\dots, \\ &{i_{l d_l}}^{n_1^{l,d_l}},\dots, {i_{l1}}^{n_1^{l,1}},
\dots, {i_{1d_1}}^{n_1^{1,d_1}},\dots, {i_{12}}^{n_1^{1,2}},
{i_{11}}^{n_1^{1,1}}),
\end{align*}
where $i^j$ means the entry $i$ appears $j$ times.  We define
$\n_j^b$ to be the subsequence
\[
\n_j^b =  (i_{b d_b}^{n_j^{b,d_b}}, \dots, i_{b1}^{n_j^{b,1}}).
\]

For a sequence $(i_1, \dots, i_m)$, let
\begin{gather*}
\ke_{i_1,\dots,i_m} = \ke_{i_1} \dots \ke_{i_m} \\
\kf_{i_1,\dots,i_m} = \kf_{i_m} \dots \kf_{i_1}.
\end{gather*}
Note the reversal of the order of the indices in the second
equation. By definition, we have $X_\infty = \ke_\n X$ and $X =
\kf_\n X_\infty$.  We define the action of $\a$ on a sequence of
vertices by $\a(i_1,\dots, i_m) = (\a(i_1), \dots, \a(i_m))$. Let
$|\n_j^b| = \sum_{c=1}^{d_b} n_j^{b,c}$ be the length of the
subsequence $\n_j^b$.  Note that
\[
|\a(\n_j^b)| = |\n_j^b|,
\]
for each $j$ and $b$.  Let $\n'$ be the sequence obtained from
$\a(\n)$ by reordering the vertices appearing in each subsequence
$\a(\n_j^b)$ such that the vertices appear in the order $i_{b
d_b}, \dots i_{b 1}$.  Note that this merely involves rearranging
vertices in the same $\a$-orbit and thus we have
\[
\ke_{\n'} = \ke_{\a(\n)},\quad \kf_{\n'} = \kf_{\a(\n)},
\]
and $|(\n')_j^b| = |\a(\n_j^b)| = |\n_j^b|$.  In particular, we
have $\ke_{\n'} X = X_\infty$.  Thus, we must have $(n')_1^{1,c}
\le n_1^{1,c}$ for $1 \le c \le d_1$ by our definition of
$n_1^{1,c}$ as the maximum number of times $\ke_{i_{1c}}$ can be
applied to $X$ without getting zero.  It follows that
$(n')_1^{1,c} = n_1^{1,c}$ for $1 \le c \le d_1$.  We can then
continue in this way to show that in fact $\n' = \n$.  Therefore,
each subsequence $\n_j^b$ is invariant under the action of $\a$
after rearranging the terms.  Thus we must have $n_j^{b,c} =
n_j^{b,c'}$ for $1 \le c,c' \le d_b$.  Denote this common value by
$n_j^b$.  It follows that
\[
\kf_\n = \kf_{\i_l}^{n_k^l} \dots \kf_{\i_1}^{n_k^1} \dots
\kf_{\i_l}^{n_1^l} \dots \kf_{\i_1}^{n_1^1}.
\]
Therefore $X = \kf_\n X_\infty \in B^g_\Gamma(\infty)$.  So $B'
\subset B^g_\Gamma(\infty)$ (as sets) which completes the proof.
\end{proof}

\begin{rem}
For the case where $Q$ is of finite type, we have by
Proposition~\ref{prop:irrcomp-finite} that the irreducible
components $B_Q^g(\infty)$ are the closures of the conormal
bundles of the various $G_\v$ orbits in $\mathbf{E}_\Omega(\v)$.
Then by a result of \cite{H04} one can see that the irreducible
components fixed by $\a$ are precisely the closures of the
conormal bundles to the orbits through the isomorphically
invariant representations (that is, the representations which are
invariant under $F(\a)$ up to isomorphism).
\end{rem}

We have realized the crystal of $U_q^-(\g(\Gamma))$ geometrically
on the set of irreducible components of Lusztig's quiver varieties
associated to $Q$ which are invariant under the automorphism $\a$.
The crystal maps are defined in terms of those coming from the
crystal of $U_q^-(\g(Q))$.  However, it is also possible to give a
more direct geometric definition of these maps.  We briefly sketch
the construction.

Recall that $\alpha_\i = \sum_{i \in \i} \alpha_i$ and that no two
vertices in the orbit $\i$ are connected by an edge. Define
\[
\mathbf{E}^0(\v)_{\i,p} = \{x \in \mathbf{E}^0(\v)\ |\
\varepsilon_i(x) = p\ \forall\ i \in \i\}.
\]
Note that $\mathbf{E}^0(\v)_{\i,p}$ is a locally closed subvariety
of $\mathbf{E}^0(\v)$.  We have for $\a$-stable $\v$
\begin{equation} \label{eq:ip-action}
F(\a)(\mathbf{E}^0(\v)_{i,p}) = \mathbf{E}^0(\v)_{\a(i),p},
\end{equation}
and therefore
\begin{equation} \label{eq:ip-def}
\mathbf{E}^0(\v)_{\i,p} = \bigcap_{i \in \i}
\mathbf{E}^0(\v)_{i,p} = \bigcap_{n=1}^{d_\i}
F(\a)^n(\mathbf{E}^0(\v)_{i,p}) \text{ for } i \in \i.
\end{equation}
Let $X \in B^g_\Gamma(\infty)$.  Thus we have $\a(X)=X$.  Then we
see from \eqref{eq:ip-action} that $\varepsilon_i(X) =
\varepsilon_j(X)$ for $i,j \in \i$ and we define their common
value by $\varepsilon_\i(X)$.  It is easy to see that this
definition coincides with the previous definition of
$\varepsilon_\i(X)$.

Now consider the diagram \eqref{eq:lus-proj1} in the case $\v =
{\bar \v} + c\alpha_\i$ for $c \in \Z_{\ge 0}$.  One easily sees
that $\mathbf{E}^0(c \alpha_\i) = \{0\}$.  Thus we have
\[
\mathbf{E}^0(\bar \v) \cong \mathbf{E}^0(\bar \v) \times
\mathbf{E}^0(c \alpha_\i) \stackrel{\varpi_1}{\longleftarrow}
\mathbf{E}^0(\bar \v,c \alpha_\i)
\stackrel{\varpi_2}{\longrightarrow} \mathbf{E}^0(\v).
\]
Using this diagram, we proceed exactly as in the case of the
quiver $Q$ and see that $B^g_\Gamma(\bar \v,\infty)_{\i,0} \cong
B^g_\Gamma(\v,\infty)_{\i,c}$ where $B^g_\Gamma(\v,\infty)_{\i,p}$
is the set of all $X \in B^g_\Gamma(\v,\infty)$ such that
$\varepsilon_\i(X)=p$.  Then we can define $\ke_\i$ and $\kf_\i$
just as we defined $\ke_i$ and $\kf_i$ and it is not hard to see
that these direct geometric definitions agree with the above
definitions in terms of compositions of the $\ke_i$ and $\kf_i$.

Let $\Lambda^\a(\v)$ be the subvariety of $\Lambda(\v)$ consisting
of the union of those irreducible components of $\Lambda(\v)$ that
are invariant under $\a$.  The above results can then be collected
into the following.

\begin{theo}
Let $B_\Gamma^g(\v,\infty)$ be the set of irreducible components
of $\Lambda^\a(\v)$ and let $B_\Gamma^g(\infty) = \bigsqcup_\v
B_\Gamma^g(\v,\infty)$.  The maps $\ke_\i$, $\kf_\i$, $\wt$,
$\varepsilon_\i$ and $\varphi_\i$, $\i \in \I$, defined above
endow $B_\Gamma^g(\infty)$ with the structure of a crystal and
this crystal is isomorphic to that associated to the crystal base
of $U_q^-(\g(\Gamma))$.
\end{theo}


\section{The crystal of an irreducible highest-weight representation
and admissible automorphisms}
\label{sec:rep-aa}

In this section and the next we apply our methods to the
irreducible highest weight representations.

Let $\w \in P(Q)^+$ be a dominant integral weight of $\g(Q)$ such
that $\a(\w) = \w$.  Thus we can also think of $\w$ as a dominant
integral weight of $\g(\Gamma)$. Recall that $B_Q(\w)$ and
$B_Q^g(\w)$ are the algebraic and geometric crystals of $V_Q(\w)$,
the irreducible highest weight representation of $U_q(\g(Q))$ of
highest weight $\w$, respectively. We know that $B_Q(\w) \cong
B_Q^g(\w)$. Let $B_\Gamma(\w)'$ denote the crystal of
$V_\Gamma(\w)$, the irreducible highest weight module of
$U_q(\g(\Gamma))$ of highest weight $\w$.  Our goal is to define a
geometric crystal $B_\Gamma^g(\w)$ of $V_\Gamma(\w)$.

As for the case of the universal enveloping algebra, the operators
$\ke_i$ and $\ke_j$ (resp. $\kf_i$ and $\kf_j$) for $i,j \in \i$
commute and we define the operators
\[
\ke_\i = \prod_{i \in \i} \ke_i,\quad \kf_\i = \prod_{i \in \i}
\kf_i,\quad \i \in \I.
\]
For $b \in B_Q(\w)$ and $\i \in \I$, we also define as before
\begin{gather*}
\varepsilon_\i (b) = \max \{k \ge 0\ |\ \ke_\i^k b \ne 0\}, \\
\varphi_\i (b) = \varepsilon_\i(b) + \left< h_\i, \wt(b)\right>.
\end{gather*}

Let $B_\Gamma(\w)$ be the subset of $B_Q(\w)$ generated by the the
$\kf_\i$, $\i \in \I$ acting on the highest weight element $b_\w$
of $B_Q(\w)$.  Note that if we restrict the map $\wt : B_Q(\w) \to
P(Q)$, where $P(Q)$ is the weight lattice of $\g(Q)$, to the
subset $B_\Gamma(\w)$, then the image lies in the subset of $P(Q)$
that is invariant under the action of $\a$.  Thus we can view it
as a map $\wt : B_\Gamma(\w) \to P(\Gamma)$ where $P(\Gamma)$ is
the weight lattice of $\g(\Gamma)$.

Note that for $b \in B_\Gamma(\w)$, $\wt(b)$ is invariant under
$\a$ and $\varepsilon_\i(b) = \varepsilon_i(b)$ for any $i \in
\i$. Thus
\begin{align*}
\varphi_\i (b) &= \varepsilon_\i(b) + \left< h_\i, \wt(b) \right>
\\
&= \varepsilon_\i(b) + \frac{1}{d_\i} \sum_{i \in \i} \left< h_i,
\wt(b) \right> \\
&= \varepsilon_i(b) + \left<h_i,\wt(b) \right> \text{ for some } i
\in \i \\
&= \max \{k \ge 0\ |\ \kf_i^k b \ne 0\} \\
&= \max \{k \ge 0\ |\ \kf_\i^k b \ne 0\}.
\end{align*}

\begin{prop}
The set $B_\Gamma(\w)$, along with the maps $\wt$,
$\varepsilon_\i$, $\varphi_\i$, $\ke_\i$ and $\kf_\i$ is a
$\g(\Gamma)$-crystal.
\end{prop}

\begin{proof}
We have to show that Axioms (1)-(7) of \cite[Definition 4.5.1]{HK}
defining a crystal are satisfied.  Axioms (1)-(6) follow easily
from the axioms for the crystal $B_Q(\w)$ and Axiom (7) is vacuous
since we never have $\varphi_\i (b) = -\infty$.
\end{proof}

We want to show that the crystal $B_\Gamma(\w)$ is actually
isomorphic to the crystal $B_\Gamma(\w)'$.

\begin{defin}
For $\w \in P(Q)^+$ (resp. $\w \in P(\Gamma)^+$), $T_\w =
\{t_\w\}$ is the crystal with one element where $\wt(t_\w) = \w$,
$\varepsilon_i(t_\w) = \varphi_i(t_\w) = -\infty$ (resp.
$\varepsilon_\i(t_\w) = \varphi_\i(t_\w) = -\infty$) and
$\ke_i(t_\w) = \kf_i(t_\w) = 0$ for all $i \in I$ (resp. $\ke_\i
(t_\w) = \kf_\i (t_\w) = 0$ for all $\i \in \I$).
\end{defin}

Note that if we view $\w \in P(\Gamma)^+$ as an $\a$-invariant
element of $P(Q)^+$, we have that the definitions of $T_\w$ for
$\g(Q)$ and $\g(\Gamma)$ coincide under the above definitions of
$\ke_\i$ and $\kf_\i$ as the composition of various $\ke_i$ and
$\kf_i$.

Let $B(\w)$ denote either $B_Q(\w)$ or $B_\Gamma(\w)'$ and
$B(\infty)$ denote either $B_Q(\infty)$ or $B_\Gamma(\infty)$. We
will use the notation $\iota$ to denote elements of either $I$ or
$\I$.  Now consider the tensor product of crystals $B(\infty)
\otimes T_\w$ and let $\pi_\w : B(\infty) \otimes T_\w \to B(\w)$
be the strict morphism given by $\kf_{\iota_1} \dots \kf_{\iota_l}
b_\infty \otimes t_\w \to \kf_{\iota_1} \dots \kf_{\iota_l} b_\w$,
where $b_\infty$ is the unique element of $B(\infty)$ of weight
zero and $b_\w$ is the unique element of $B(\w)$ of weight $\w$.
We have that
\[
\pi_\w : \{b \in B(\infty) \otimes T_\w\ |\ \pi_\w(b) \ne 0\}
\cong B(\w).
\]

The following characterization of the crystal of an irreducible
highest weight representation with be useful.

\begin{prop}[\cite{KS97}] \label{prop:crystal-irrep-char}
Let $B$ be a $\g(Q)$ or $\g(\Gamma)$-crystal and $b_\w$ an element
of $B$ with weight $\w \in P^+$.  Assume the following conditions.
\begin{enumerate}
\item \label{irrep-char1} $b_\w$ is the unique element of $B$ with
weight $\w$.

\item \label{irrep-char2} There is a strict morphism $\Phi :
B(\infty) \otimes T_\w \to B$ such that $\Phi (b_\infty \otimes
t_\w) = b_\w$.

\item \label{irrep-char3} The set $\{b \in B(\infty) \otimes T_\w\
|\ \Phi(b) \ne 0 \}$ is isomorphic to $B$ through $\Phi$ as a set.

\item \label{irrep-char4} For any $b \in B$, $\varepsilon_\iota(b)
= \max \{k \ge 0\ |\ \ke_\iota^k(b) \ne 0\}$ and $\varphi_\iota(b)
= \max \{k \ge 0\ |\ \kf_\iota^k(b) \ne 0\}$ for all $\iota$.
\end{enumerate}
Then $B$ is isomorphic to $B(\w)$.
\end{prop}

\begin{prop}
\label{prop:rep-crystal}
We have $B_\Gamma(\w) \cong
B_\Gamma(\w)'$.
\end{prop}

\begin{proof}
We prove this by an appeal to
Proposition~\ref{prop:crystal-irrep-char}, whose conditions are
easily verified.
\end{proof}


\section{A geometric realization of the crystal of an irreducible
highest-weight representation in the non-simply laced case}
\label{sec:geom-real-nak}

Since we know that the geometric crystal $B^g_Q(\w)$ is isomorphic
to the algebraic crystal $B_Q(\w)$, we can define $B^g_\Gamma(\w)$
just as we defined $B_\Gamma(\w)$ and we know that
$B^g_\Gamma(\w)$ is isomorphic to the crystal associated to the
crystal base of $V_\Gamma(\w)$.  By definition, the underlying set
of $B^g_\Gamma(\w)$ consists of the irreducible components of
Lusztig's quiver varieties generated from the highest weight
component (which consists of a point) by the action of the
operators $\kf_\i$ for $\i \in I$.  We now seek to describe this
set in a more direct and geometric way.

We first extend the action of the automorphism $\a$ to Nakajima's
quiver varieties as follows.  Let $\V$ and $\W$ be $I$-graded
vector spaces of graded dimension $\v$ and $\w$ respectively.  As
before, we define the new vector spaces $^\a \V$ and $^\a \W$ and
we have the action $F(\a) : \Lambda(\v) \to \Lambda(\a(\v))$ on
Lusztig's quiver varieties.  For an element $t=(t_i) \in
\bigoplus_{i \in I} \Hom (\V_i,\W_i)$, we define a new element
${^\a t} = (^\a t_i) \in \bigoplus_{i \in I} \Hom ({^\a \V}_i,
{^\a \W}_i)$ by ${^\a t}_i = t_{\a^{-1}(i)}$.  Combining these
action yields a map
\[
F(\a) : \Lambda(\v,\w) \to \Lambda(\a(\v),\a(\w)),
\]
which in turn induces a map
\[
F(\a) : \L(\v,\w) \to \L(\a(\v),\a(\w)).
\]
Note that if we restrict ourselves to $V$ and $W$ whose dimensions
are invariant under $\a$, we can consider $F(\a)$ as an
automorphism of Nakajima's quiver variety and thus an automorphism
of the set of irreducible components.  As before, we denote the
automorphism of the set of irreducible components thus obtained by
$\a$.

\begin{prop}
The underlying set of $B^g_\Gamma(\w)$ is precisely the subset of
irreducible components of $B^g_Q(\w)$ consisting of those
components that are invariant under the automorphism $\a$.
\end{prop}

\begin{proof}
The proof is almost exactly analogous to that of
Proposition~\ref{prop:univ-irrcomp-inv} and will be omitted.
\end{proof}

As in the case of Lusztig's quiver variety, we should note that
the irreducible components invariant under $\a$ are not
necessarily invariant under $F(\a)$. That is, each individual
point in the irreducible component is not necessarily invariant
under $F(\a)$.

We have realized the crystal of $V_\Gamma(\w)$ geometrically on
the set of irreducible components of the Lusztig quiver variety
associated to $Q$ which are invariant under the automorphism $\a$.
The crystal maps are defined in terms of those coming from the
crystal of $V_Q(\w)$.  However, as for the case of the universal
enveloping algebra, it is also possible to give a more direct
geometric definition of these maps. We briefly sketch the
construction.

Recall that $\alpha_\i = \sum_{i \in \i} \alpha_i$ and that no two
vertices in the orbit $\i$ are connected by an edge. Define
\[
\L(\v,\w)_{\i,p} = \{[x,t] \in L(\v,\w)\ |\ \varepsilon_i((x,t)) =
p\ \forall\ i \in \i\}.
\]
Note that $\L(\v,\w)_{\i,p}$ is a locally closed subvariety of
$\L(\v,\w)$.  We have for $\a$-stable $\v$ and $\w$
\begin{equation} \label{eq:ip-action-nak}
F(\a)(\L(\v,\w)_{\i,p}) = \L(\v,\w)_{\a(i),p},
\end{equation}
and therefore
\begin{equation} \label{eq:ip-def-nak}
\L(\v,\w)_{\i,p} = \bigcap_{i \in \i} \L(\v,\w)_{i,p} =
\bigcap_{n=1}^{d_\i} F(\a)^n(\L(\v,\w)_{i,p}) \text{ for } i \in
\i.
\end{equation}
Let $X \in B^g_\Gamma(\w)$.  Thus we have $\a(X)=X$.  Then we see
from \eqref{eq:ip-action-nak} that $\varepsilon_i(X) =
\varepsilon_j(X)$ for $i,j \in \i$ and we define their common
value by $\varepsilon_\i(X)$.  It is easy to see that this
definition coincides with the previous definition of
$\varepsilon_\i(X)$.

Now consider the diagram \eqref{eq:diag_action_mod} in the case
$\v = \v' + c\alpha_\i$ for $c \in \Z_{\ge 0}$.  Then we have
\[
\L(\v',\w) \stackrel{\pi_1}{\longleftarrow} \F(\v,\w,c \alpha_\i)
\stackrel{\pi_2}{\longrightarrow} \L(\v,\w).
\]
Using this diagram, we proceed exactly as in the case of the
quiver $Q$ and see that $B^g_\Gamma(\v',\w)_{\i,0} \cong
B^g_\Gamma(\v,\w)_{\i,c}$ where $B^g_\Gamma(\v,\w)_{\i,p}$ is the
set of all $X \in B^g_\Gamma(\v,\w)$ such that
$\varepsilon_\i(X)=p$.  Then we can define $\ke_\i$ and $\kf_\i$
just as we defined $\ke_i$ and $\kf_i$ and it is not hard to see
that these direct geometric definitions agree with the above
definitions in terms of compositions of the $\ke_i$ and $\kf_i$.

Let $\L^\a(\v,\w)$ be the subvariety of $\L(\v,\w)$ consisting of
the union of those irreducible components of $\L(\v,\w)$ that are
invariant under $\a$.  The above results can then be collected
into the following.

\begin{theo}
Let $B_\Gamma^g(\v,\w)$ be the set of irreducible components of
$\L^\a(\v,\w)$ and let $B_\Gamma^g(\w) = \bigsqcup_\v
B_\Gamma^g(\v,\w)$.  The maps $\ke_\i$, $\kf_\i$, $\wt$,
$\varepsilon_\i$ and $\varphi_\i$, $\i \in \I$, defined above
endow $B_\Gamma^g(\w)$ with the structure of a crystal and this
crystal is isomorphic to that associated to the crystal base of
$V_\Gamma(\w)$, the irreducible highest weight module of
$U_q(\g(\Gamma))$ of highest weight $\w$.
\end{theo}


\section{A geometric realization of the spin representation of
$\mathfrak{so}_{2n+1}$} \label{sec:spinrep}

In this section, we will apply the results of this paper to give
an explicit realization of the spin representation of the Lie
algebra $\mathfrak{so}_{2n+1} = \mathfrak{so}_{2n+1}\C$ of type
$B_n$. In the process, we will end up with a combinatorial
description of the crystal in terms of Young diagrams.

\begin{figure}
\centering \epsfig{file=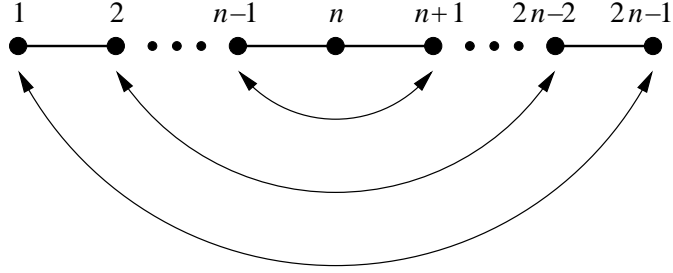,width=0.7\textwidth}
\caption{The admissible automorphism of the Dynkin diagram of type
$A_{2n-1}$ realizing the valued graph of type $B_n$.  The arrows
indicate the action of the automorphism. \label{fig:quiveraut}}
\end{figure}

Let $\g = \mathfrak{so}_{2n+1}$.  The spin representation is the
highest weight representation $V(\omega_n)$ of highest weight
$\omega_n$.  The valued graph, or Dynkin diagram, of type $B_n$
can be obtained from the Dynkin diagram of type $A_{2n-1}$ and the
automorphism $\a$ as depicted in Figure~\ref{fig:quiveraut}.
We let $Q$ be the quiver of type $A_{2n-1}$ (we pick an
orientation compatible with the automorphism) and $\Gamma$ the
corresponding valued graph of type $B_n$. Since we are interested
in the spin representation, we should consider $\w=\e^n$ where
$\e^n \in (\Z_{\ge 0})^{2n-1}$ is the element with $n$th component
equal to one and all other components equal to zero. This element
is invariant under the automorphism $\a$.

\begin{figure}
\centering \epsfig{file=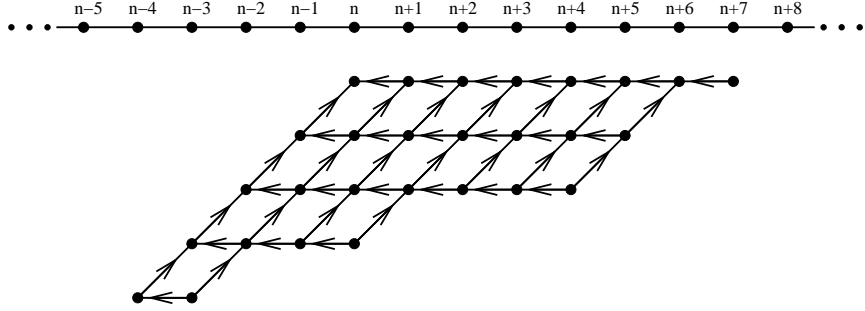,width=0.9\textwidth}
\caption{The irreducible component $X_Y$ corresponding to the
Young diagram $(8,7,7,4,2)$.  The Dynkin diagram is shown above
the Young diagram.  Each vertex in the Young diagram corresponds
to a basis vector in the representation of the quiver.  Its degree
is given by the vertex of the Dynkin diagram above it.  The arrows
in the Young diagram indicate the action of $x$ (i.e. $x$ takes
the vector at the tail of an arrow to the vector at the tip and
all other components of $x$ are zero). \label{fig:yd-irrcomp}}
\end{figure}

Let $\mathcal{Y}$ be the set of all Young diagrams (or partitions)
$Y = (n \ge \lambda_1 \ge \lambda_2 \ge \dots \ge \lambda_n)$ with
at most $n$ rows (or parts) and with all rows of length at most
$n$.  That is, $\mathcal{Y}$ is the set of all Young diagrams that
fit inside an $n \times n$ box. To a Young diagram $Y \in
\mathcal{Y}$, we can associate an irreducible component $X_Y$ as
in \cite{FS03}.  Such an irreducible component can be pictured as
in Figure~\ref{fig:yd-irrcomp}.

\begin{prop}
The association $Y \leftrightarrow X_Y$ is a one-to-one
correspondence between the set $\mathcal{Y}$ and the set of
irreducible components $B^g_Q(\w)$.  Furthermore, each irreducible
component is a point.
\end{prop}
\begin{proof}
The case of type $A_\infty$ is done in \cite[Thm~5.1]{FS03}. The
case of type $A_{2n-1}$ is exactly analogous (see
\cite[Rem~5.4]{FS03}).  The fact that the irreducible components
are points can either be seen directly or from the dimension
formula for quiver varieties (see \cite[Cor~3.12]{N98}).
\end{proof}

For a Young diagram $Y$, let $\tilde Y$ denote the conjugate Young
diagram.  Thus if $Y = (\lambda_1 \ge \dots \ge \lambda_m)$ then
${\tilde Y} = (\mu_1 \ge \dots \ge \mu_k)$ where $\mu_i = \#\{j\,
|\, \lambda_j \ge i\}$.  Pictorially, $\tilde Y$ is obtained from
$Y$ by reflection in the diagonal.

\begin{prop}
We have $\a(X_Y) = X_{\tilde Y}$.
\end{prop}
\begin{proof}
This is easy to see from our explicit description of the
irreducible components $X_Y$.  It also follows from dimension
arguments since there is at most one irreducible component of a
given dimension.
\end{proof}

Let $\mathcal{Y}^\a$ be the set of self-conjugate elements of
$\mathcal{Y}$.  Then we have that
\[
B^g_\Gamma(\w) = \{X_Y\, |\, Y \in \mathcal{Y}^\a\}.
\]
Now, as can be seen in Figure~\ref{fig:yd-irrcomp}, the Young
diagrams in $\mathcal{Y}$ can be thought of as consisting of
vertices of various degrees (between 1 and $2n-1$).  For a Young
diagram $Y \in \mathcal{Y}$, let $Y^+_k$ (resp. $Y^-_k$) be the
Young diagram in $\mathcal{Y}$ obtained from $Y$ by adding (resp.
removing) a vertex of degree $k$ if such a Young diagram exists
and let $Y^\pm_k$ be the ghost partition $\nabla$ if such a Young
diagram does not exist. Using the fact that the quiver varieties
are points, it is easy to see (see \cite{FS03,S03}) that the
action of the Kashiwara operators $\ke_k$ and $\kf_k$ on
$B^g_Q(\w)$ is given by $\ke_k(X_Y) = X_{Y^-_k}$ and $\kf_k(X_Y) =
X_{Y^+_k}$ where $X_\nabla = 0$.

For $Y \in \mathcal{Y}$, let $Y^{++}_{i,j}$ (resp. $Y^{--}_{i,j}$)
be the Young diagram in $\mathcal{Y}$ obtained from $Y$ by adding
(resp. removing) a vertex of degree $i$ and a vertex of degree
$j$.  Then we have the following theorem.

\begin{theo}
The set $\mathcal{Y}^\a$ with operators
\begin{gather*}
\ke_k (Y) = Y^{--}_{k,2n-k},\quad \kf_k = Y^{++}_{k,2n-k},\quad 1
\le
k \le n-1, \\
\ke_n(Y) = Y^-_n,\quad \kf_n(Y) = Y^+_n, \\
\wt (Y) = \omega_n - \sum_{i=1}^n \#\{\text{vertices in $Y$ of
degree $i$}\} \alpha_i, \\
\varepsilon_k(Y) = \max \{k\, |, \ke_k(Y) \ne 0\}, \quad
\varphi_k(Y) = \max \{k\, |\, \kf_k(Y) \ne 0\},
\end{gather*}
where we identify the ghost partition $\nabla$ with zero, is
isomorphic to the crystal of the spin representation of type
$B_n$.
\end{theo}
\begin{proof}
This follows immediately from Proposition~\ref{prop:rep-crystal}
and the above discussion.  We have identified $Y \in
\mathcal{Y}^\a$ with the irreducible component $X_Y$.
\end{proof}

\begin{prop}
Let $V$ be the vector space spanned by $\mathcal{Y}^\a$, let $e_k$
and $f_k$ act on this space by extending the action of $\ke_k$ and
$\kf_k$ by linearity, and define $h_k(Y) = \left< h_k, \wt(Y)
\right> Y$ (extended by linearity).  Then $V \cong V(\omega_n)$ as
representations of $\mathfrak{so}_{2n+1}$ where
$\{e_k,f_k,h_k\}_{k=1}^n$ are the usual Chevalley generators. That
is, we can realize the spin representation in a natural way on the
space spanned by the self-conjugate Young diagrams fitting inside
an $n \times n$ box.
\end{prop}
\begin{proof}
This can be shown by direct computation.  One merely checks that
the necessary commutation relations hold.
\end{proof}

This realization of the spin representation of type $B_n$ is
similar to the realization of the spin representation of type
$D_n$ on a certain set of Young diagrams obtained in \cite{S03}.


\bibliographystyle{abbrv}
\bibliography{biblist}

\end{document}